\documentclass[11pt,a4paper,twoside]{article}  
\usepackage[cp1251]{inputenc}                            
\usepackage{amsmath,amsfonts,amssymb,amscd,euscript}
\usepackage[english, russian]{babel}

\usepackage{graphicx}
\usepackage{psfrag}
\usepackage{cite}
\makeatletter

\def\titlerus{\thispagestyle{empty} { } \vspace{-5mm} \noindent
\raisebox{-37pt}[\headheight][0pt]{\vbox{ \hbox to \textwidth{\hfil
\scriptsize 
\ \hfill \
}
\vspace{2pt} \hrule \vspace{8pt} \hbox to \textwidth{\series \hfil  \issue}
\vspace{30pt} \hbox{УДК \UDK} }} \vspace{ 30pt plus 6pt }}

\def\titleeng{\vspace{3ex} \hfill Поступила  в редакцию  \ \datereceive \par \vspace{5ex} \par
\noindent \parbox{166mm
}{\small \noindent \contactinformation
\par \vspace{30pt} \par {\textbf {\textit {\autorseng}}} \par {\bf \articleseng} \par
\vspace{10pt}\par {\it Keywords}: \keywordseng
\par \vspace{7pt} \par \noindent \small {Mathematical Subject Classifications}: \MSC }
\par \vspace{20pt}  \par \small  \noindent \referateng}

\def\received{\vspace{3ex} \hfill Received \ \datereceive \par \vspace{5ex} }

\def\annotationandkeywordsrus{\noindent {\small \referatrus \par } \vspace{8pt}
\noindent {\small {\it Ключевые слова}: \keywordsrus} \par \vspace{10pt}}

\@addtoreset{equation}{section}
\renewcommand{\section}{\@startsection{section}{1}{0pt}{1.3ex
plus 1ex minus .1ex}{1.3ex plus .1ex}{\bf\,\S\,}}
\newcommand{\point}{\hspace*{-4mm}{\bf.}\;}
\newcommand{\sect}[1]{\begin{flushleft}%
\protect{\section{\point#1}}\end{flushleft}}

\renewcommand{\@begintheorem}[2]{\begin{trivlist}
\item[\hspace{\labelsep}{\bf \mbox{~~~}#1\ #2.}]}
\renewcommand{\@opargbegintheorem}[3]{\begin{trivlist}
\item[\hspace{\labelsep}{\bf \mbox{~~~}#1\ #2 {\rm (#3).}}]}
\renewcommand{\@endtheorem}{\end{trivlist}}

\newtheorem{teo}{Теорема}

\newtheorem{pre}{Предложение}

\newtheorem{sle}{Следствие}

\newtheorem{zam}{Замечание}

\newcommand{\doc}{\mbox{Д о к а з а т е л ь с т в о}}
\renewcommand{\@evenfoot}{}
\renewcommand{\@oddfoot}{}

\renewcommand*{\@biblabel}[1]{#1.\hfill}

\newcommand*{\CSep}{.\ }
\renewcommand{\@makecaption}[2]{%
  \vskip\abovecaptionskip
  \sbox\@tempboxa{{\bf #1\CSep}{#2}}%
  \ifdim \wd\@tempboxa >\hsize
  \begin{center}%
    {\footnotesize{\bf #1\CSep}{#2\par}}%
  \end{center}%
  \else
    \global \@minipagefalse
    \hb@xt@\hsize{\hfil\box\@tempboxa\hfil}%
  \fi
  \vskip\belowcaptionskip%
}

\renewcommand{\@evenhead}{\raisebox{0pt}[\headheight][0pt]{\vbox{\hbox to\textwidth{\thepage \strut \hfil
\text{\autorsrus} \hfil } \hrule \vspace{8pt} \hbox to \textwidth{\series \hfil  \issue}}}}

\renewcommand{\@oddhead}{\raisebox{0pt}[\headheight][0pt]{\vbox{\hbox to\textwidth{  \strut \hfil
\text{\articleshortname} \hfil \thepage} \hrule \vspace{8pt} \hbox to \textwidth{\series \hfil  \issue}}}}

\newcommand{\series}{
\ }

\newcommand{\issue}{2013. Вып.\,} 

\newcommand{\autorsrus}{Д.\,В.~Хлопин} 
\newcommand{\autorseng}{D.\,V.~Khlopin} 

\newcommand{\articleshortname}{Краевые  условия на бесконечности для
строгого оптимального управления }
\newcommand{\articleseng}{On necessary boundary conditions for strictly optimal control in infinite horizon control problems}

\newcommand{\UDK}{\ 517.977.5
} 
\newcommand{\MSC}{49K15,49J45,37N40} 
\newcommand{\referatrus}{
В работе   рассматриваются задачи управления на бесконечном промежутке со
свободным правым концом. Получены необходимые условия строгой оптимальности. Сам метод доказательства фактически следует классической работе Халкина, а построенное в работе краевое условие на бесконечности является усилением  условия, предложенного Сейерстадом.
   Построенная в работе полная система соотношений принципа максимума позволяет выписать для
   сопряженной переменной выражение в виде несобственного интеграла, зависящего лишь от
   разворачивающейся траектории. С.М.Асеев, А.В.Кряжимский,
   V.M.Veliov получали такое выражение   в качестве необходимого в некоторых классах задач управления.
   Строгая оптимальность в ряде случаев позволяет создать переопределенную систему соотношений, в работе получены условия, достаточные для этого. Разобран пример.
}

\newcommand{\referateng}{
In the paper we consider the infinite horizon control problems on the interval with free right-hand endpoint. We obtain the necessary conditions of strict optimality. The method of the proof actually follows the classic paper by Halkin, and the boundary condition for infinity that we construct in our paper is a stronger variety of the Seierstad condition. The complete system of relations of the maximum principle that was obtained in the paper allows us to write the expression for the adjoint variable in the form of improper integral that depends only on the developing trajectory. S.M. Aseev, A.V. Kryazhimskii, and V.M. Veliov obtained the similar condition as a necessary condition for certain classes of control problems. As we note in our paper, the obtained conditions of strict optimality lead us to a redefined system of relations for sufficiently broad class of control problems. An example is considered.
   }

\newcommand{\keywordsrus}{задача управления, строго оптимальное управление;
 задача на бесконечном промежутке,  необходимые условия оптимальности,  краевое условие
на бесконечности, принцип максимума Понтрягина
  }
\newcommand{\keywordseng}{ Optimal control; infinite horizon problem, transversality condition
for infinity, necessary conditions of optimality, strong optimal control;
 Pontryagin Maximum Principle }

\newcommand{\datereceive}{26.12.12} 

\newcommand{\contactinformation}{Хлопин Дмитрий Валерьевич,
к.\,ф.-м.\,н., Институт математики и механики УрО РАН, 620990,
Россия, г. Екатеринбург, ул. Софьи Ковалевской,\,16, E-mail:
khlopin@imm.uran.ru}
\newcommand{\contactinformationenglish}{Khlopin Dmitry Valer'evich,
Candidate of Physics and Mathematics, Senior Researcher,
    Institute of Mathematics and Mechanics, Ural Branch,     Russian
Academy of Sciences, 16, S.Kovalevskaja St., 620990, Yekaterinburg, Russia. \\
E-mail: khlopin@imm.uran.ru }

\DeclareMathOperator{\esslim}{ess}
\usepackage[active]{srcltx}

\newcommand{\rav}{\stackrel{\triangle}{=}}
\newcommand{\rref}[1]{$(\ref{#1})$}
\newcommand{\mm}[1]{{\mathbb{#1}}}
\newcommand{\ct}[1]{{\mathcal{#1}}}
\newcommand{\fr}[1]{{\mathfrak{#1}}}
\newcommand{\bo}{
{\hfill {$\Box$}} }
\begin{document}
\sloppy

\titlerus

\begin{flushleft}
{\bf \copyright { \textit { \autorsrus}} \\[2ex]
{ О НЕОБХОДИМЫХ КРАЕВЫХ УСЛОВИЯХ ДЛЯ СТРОГОГО ОПТИМАЛЬНОГО УПРАВЛЕНИЯ
 В ЗАДАЧАХ УПРАВЛЕНИЯ НА БЕСКОНЕЧНОМ ПРОМЕЖУТКЕ
}
\footnote{Работа частично поддержана  грантом  №~12-01-31172
}}
\end{flushleft}

\annotationandkeywordsrus

\sloppy
  \sect{Определения и обозначения}

  Рассмотрим задачу  максимизации на бесконечном промежутке функционала
   \begin{equation}
   \label{opt}
   J[u](T)\rav\int_{[0,T]} g(t,x(t),u(t))\, dt\ \text{ при } \
   T\to\infty.
   \end{equation}
   на траекториях управляемой системы
  \begin{subequations}
\begin{eqnarray}
\label{sys1}    \dot{x}=f(t,x,u)&\ &\text{для п.в.}\ t\in\mm{T},\\
\label{sys2}     x(0)=x_*,&\ & \\
\label{sys4}     u\in U(t)&\ &\text{для п.в.}\ t\in\mm{T}.
   \end{eqnarray}
  \end{subequations}
   Здесь в качестве промежутка времени рассматривается полуось
  $\mm{T}\rav
  \mm{R}_{\geqslant 0};$
     в качестве фазового пространства управляемой системы используется некоторое
   конечномерное
 евклидовое пространство $\mm{X}$.
   Также считается заданным конечномерное евклидовое пространство $\mm{U}$ и
   многозначное   отображение  $U:\mm{T}\rightsquigarrow\mm{U}$.
   Под множеством всевозможных допустимых управлений $\fr{U}$
   понимается
  множество всех измеримых по Борелю локально ограниченных селекторов многозначного отображения $U$.

  Будем всюду далее предполагать, что $U$ компактнозначно,
  локально интегрально ограничено с измеримым по Борелю графиком;
   $f,g$ и их производные по $x$~---
  локально липшицевые по $x$ отображения Каратеодори,  интегрально
  ограниченные на компактах; кроме того, для $f$ выполнено условие подлинейного роста:
  для некоторой непрерывной функции $M\in C(\mm{T},\mm{T})$
  $$||f(t,x,u)||\leqslant M(1+||x||)\quad \forall t\in\mm{T},u\in U(t), x\in\mm{X}.$$

  Отметим, что вообще говоря,
  интеграл в \rref{opt} не обязан ни при всех допустимых управлениях игроков, ни при каких-то их управлениях,
  при $T\to\infty$ сходиться ни к
  конечному числу, ни к бесконечности. В этом случае требуется
  формализовать, что понимать под оптимальностью на бесконечности.
  Первые такие формализации появились уже  в модели Рамсея; см. \cite{4,10}.

  В задачах управления на бесконечном промежутке предложена целая
  серия таких определений оптимальности \cite{slovak,car1,car,stern}.
  Нас будет интересовать понятие строго оптимального управления
  ("strongly optimal"; см., например, \cite{car1,car}).

 Будем говорить, что оптимальное управление  $u^0\in\fr{U}$
 является {\it строго оптимальным}, если
  для всякого момента $n\in\mm{N}$ найдется не меньший его момент времени $\tau_n,$
  для которого $J[u^0](\tau_n)\geqslant J[\overline{u}](\tau_n)$
 при всех $\overline{u}\in\fr{U}$.

 Всюду далее будем предполагать, что
 строго оптимальное управление  $u^0\in\fr{U}$ существует;
 соответствующую ему траекторию обозначим через $x^0$.
 Отметим, что в самом определении уже
 построена неограниченно возрастающая
 последовательность
 моментов времени $\tau_n$, вдоль которой собственно  и происходит
 оптимизация; всюду далее будем, при необходимости, говорить о
 $\tau$-строгом оптимальном управлении.

   Отметим, что вообще говоря в работах \cite{car1,car} при определении строго оптимального решения
   накладывалось также
    дополнительное требование ограниченности функционала $J$.

\sect{Соотношения принципа максимума Понтрягина}

 Определим функцию Гамильтона--Понтрягина
  $\ct{H}:\mm{X}\times \mm{T}\times {U}_{all}\times\mm{T}\times\mm{X}\mapsto\mm{R}$
  правилом:
   $\ct{H}(x,t,u,\lambda,\psi)\rav\psi f\big(t,x,u\big)+\lambda
   g\big(t,x,u\big).$
 Введём соотношения
\begin{subequations}
 \begin{eqnarray}
   \label{sys_psi}
       \dot{\psi}(t)&=&-\frac{\partial
       \ct{H}}{\partial x}\big(x^0(t),t,u^0(t),\lambda,\psi(t)\big);\\
   \label{maxH}
\sup_{p\in U(t)}\ct{H}\big(x^0(t),t,p,\lambda,\psi(t)\big)&=&
        \ct{H}\big(x^0(t),t,u^0(t),\lambda,\psi(t)\big);\\
   \label{dob}
   \lambda=1\qquad\text{или}& & \lambda=0,\ ||\psi(0)||=1.
   \end{eqnarray}
\end{subequations}
  Хорошо известно  \cite{Halkin}, что эти соотношения выполнены для
  пары $(x^0,u^0)$ при некоторых $\psi^0\in
  (AC)(\mm{T},\mm{X}),\lambda^0\in[0,1]$. При этом от $u^0$
  требуется лишь неулучшаемость данного управления на некоторых, сколь угодно больших отрезках $[0,T]$
  при дополнительном фиксировании правого конца $x^0(T),$ что для строго оптимального управления заведомо выполнено.

   Однако, как отмечено уже в той же статье \cite{Halkin},
   построенная система соотношений \rref{sys_psi}--\rref{dob} неполна, и, вообще говоря, требует
   дополнительного условия на бесконечности. Найти сколько-нибудь
    универсальное условие трансверсальности
  хотя бы  для задачи со свободным правым концом пока не
   удается:
 большинство предложенных условий могут
выделить слишком много решений  соотношений
\rref{sys_psi}--\rref{dob}, или
 оказаться несовместными с ними (см. \cite[\S~6 и \S~16]{kr_as},\,\cite{norv}).

     С другой стороны, в работе \cite{norv} для достаточно широкого класса
     задач управления при решении соотношений принципа максимума A. Seierstad предложил
      выделять то
     из них,
     что является пределом решений этих же соотношений, но выписанных
     для задач,
     укороченных во все более поздние моменты времени.
     В самой работе  \cite{norv} необходимость такого условия
     \cite[Theorem~8.1]{norv}
     показана для весьма широкого класса задач управления, но при достаточно сильных
     предположениях,
     обеспечивающих в частности,
     суммируемость целевого функционала, а также
     сходимость $\psi^0$ к нулю
     (подробнее см., например, \cite[\S12]{kr_as}).

    В рамках этого подхода для задачи на бесконечном промежутке со свободным правым концом
      стоит искать решение
     $(\psi^0,\lambda^0)$
     соотношений принципа максимума  \rref{sys_psi}--\rref{dob}, у которого сопряженная переменная
     $\psi^0$
     является поточечным пределом
     сопряженных переменных $\psi^n$~--- решений \rref{sys_psi}~--- но зануляющихся вдоль неограниченно возрастающей последовательности
      моментов времени.

     Определим подобное асимптотическое условие для строго
     оптимального решения.
     В отличие от работы \cite{norv} мы не будем предполагать (или доказывать)
     нормальность задачи, поэтому в анормальных задачах переходить к пределу придется не
     только в компоненте $\psi$, но и в компоненте $\lambda.$
     С другой стороны, в отличие от \cite{norv} мы сможем
     гарантировать, что пары $(\lambda^n,\psi^n)$ удовлетворяют не
     только
     сопряженной системе \rref{sys_psi}, но и условию максимума
      \rref{maxH}. Формализуем сказанное.

   Будем говорить, что нетривиальное решение
   $(\lambda^0,\psi^0)$ соотношений принципа максимума \rref{sys_psi}--\rref{maxH}
  {\it строго $\tau$-исчезающее} (или просто
   строго исчезающее),
   если
    для
 некоторой подпоследовательности $\tau'$ последовательности $\tau$
существует сходящаяся к $(\lambda^0,\psi^0)$ равномерно на всяком
компакте
   последовательность пар
   $(\lambda^n,\psi^n)\in[0,1]\times (AC)(\mm{T},\mm{X})$,
   удовлетворяющих каждая
   соотношениям \rref{sys_psi}, \rref{maxH} на $[0,\tau'_n]$,  а
   также
   краевому условию
 \begin{equation}
   \label{dob_k}
   \psi^n(\tau'_n)=0.
 \end{equation}

 Такое определение названо строго исчезающим решением
 соотношений принципа максимума по аналогии с показанным в работе
 \cite{mymy} для критерия слабо догоняющее оптимальное управление
 ("weakly  uniformly overtaking optimal") более слабым необходимым условием оптимальности,
 названным там "исчезающим решением".

     Покажем при сделанных выше предположениях
     необходимость подобного условия для всякого строго оптимального решения. Само доказательство при этом фактически следует работе \cite{Halkin}, лишь дополнительно отслеживается краевое условие.


\begin{teo}
\label{1}{\it
    Пусть   $u^0\in\fr{U}$~--- $\tau$-строгое оптимальное управление в задаче
    \rref{opt}--\rref{sys4}.
    Тогда
    найдётся нетривиальное строго $\tau$-исчезающее решение
   $(\lambda^0,\psi^0) \in\{0,1\}\times (AC)(\mm{T},\mm{X})$
    соотношений принципа максимума \rref{sys_psi}--\rref{maxH},
    то есть\\
   существует сходящаяся к нему равномерно на всяком компакте
   последовательность пар
   $(\lambda^n,\psi^n)\in[0,1]\times (AC)(\mm{T},\mm{X})$
     и такая подпоследовательность $\tau'\subset\tau$, что
     каждая пара  $(\lambda^n,\psi^n)$ удовлетворяет на $[0,\tau'_n]$ соотношениям принципа максимума \rref{sys_psi}--\rref{dob}
   с краевым условием \rref{dob_k}.

   Если при этом $\lambda^0=1$, то можно считать, что все
   $\lambda^n=1$.
}
\end{teo}
\doc.\
   Зафиксируем $n\in\mm{N}.$
   Из определения строгой оптимальности   следует, что
   $J[u^0](\tau'_n)\geqslant J[u](\tau'_n)$ для всех допустимых
   $u\in\fr{U}$.
   Тогда $u_0|_{[0,\tau_n]}$ является  оптимальным при
   максимизации
   \begin{equation*}
   \label{opt_}
   J[u](\tau_n)\rav\int_{[0,\tau_n]} g(t,x(t),u(t))\, dt\
   \end{equation*}
     на траекториях \rref{sys1}--\rref{sys4}.
    По  принципу максимума  \cite{ppp,clarke}, для этой задачи найдётся
    пара
   $(\lambda^n,\psi^n)\in[0,1]\times (AC)([0,\tau_n],\mm{X}),$
   для которой на промежутке $[0,\tau_n]$ выполнены соотношения
    \rref{sys_psi},
          \rref{maxH}, а кроме того выполнено условие транcверсальности на правый конец, то есть
          $\psi^n(\tau_n)=0$.
   Более того, нормируя при необходимости, можно считать, что
   \begin{equation}
   ||\psi^n(0)||+\lambda^n=1.
      \label{dob_kk}
   \end{equation}
   Продолжим в силу управления $u^0$ функцию $\psi^n$ как решение
   сопряженной системы \rref{sys_psi} на всю ось $\mm{T}.$

   Заметим, что в силу условия \rref{dob_kk}
   $(\psi^n(0),\lambda^n)$ лежат на компакте, следовательно, перейдя при необходимости к подпоследовательности,
   можно считать, что  последовательность
   пар
   $(\psi^n(0),\lambda^n)$ сходится
   к некоторой $(\overline{\psi},\overline{\lambda})$. Рассмотрим решение $\psi^\infty$ сопряженной
    системы
    \rref{sys_psi} c
   начальным условием $\psi^\infty(0)=\overline{\psi}$ при $\lambda^\infty=\overline{\lambda}.$
   По теореме о непрерывной зависимости решений дифференциальных уравнений от  начальных данных
   имеем, что  последовательность
    $\psi^n$ сходится к $\psi^\infty$ равномерно на каждом компакте.
    Для почти каждого $t\in\mm{T}$ условие максимума \rref{maxH} имеет место для всех $(\lambda^n,\psi^n)$,
    начиная с некоторого $n$. Переходя к пределу, имеем условие максимума
    \rref{maxH} уже и для $(\lambda^\infty,\psi^\infty)$ при почти
    всех $t\in\mm{T}$.
     Построенное решение, в силу предельного перехода в \rref{dob_kk},
     удовлетворяет
     $||\psi^\infty(0)||+\lambda^\infty=1,$ в частности нетривиально.

Построим теперь строго $\tau$-исчезающее решение
$(\psi^0,\lambda^0)$, удовлетворяющее \rref{dob}. Действительно, если
$\lambda^\infty=1$, то примем
   $\psi^0\rav\psi^\infty,\lambda^0\rav 1.$
Если $\lambda^\infty=0$, то
   дополнительно поделим для каждого $n\in\mm{N}$ функции
    $\psi^n,\lambda^n$ на $||\psi^n(0)||;$
   тогда они сойдутся к
   $\psi^0\rav\psi^\infty/||\psi^\infty(0)||,\lambda^0\rav 0$.
   Если же  $\lambda^\infty\neq 0$, то
      дополнительно поделим для каждого $n\in\mm{N}$ функции $\psi^n,\lambda^n$ на
       $\lambda^n;$ тогда они сойдутся к
    $\psi^0\rav\psi^\infty/(1-||\psi^\infty(0)||),\lambda^0\rav 1.$
  Во всех случаях условие \rref{dob} выполнено. \bo

 Покажем, что можно найти дополнительное условие уже на саму
 последовательность $\tau.$

 \begin{sle}
\label{aff1} {\it
    Для всякого $\tau$-строго оптимального управления $u^0\in\fr{U}$
    и порожденной им траектории $x^0\in (AC)(\mm{T},\mm{X}),$
      для всех $n\in\mm{N}$, кроме, быть может, конечного их числа,
      выполнено также
 \begin{equation}
          \label{maxH_k}
\esslim\,\lim_{t\to \tau_n-0} \sup_{p\in
U(t)}\big[g\big(t,x^0(t),p\big)-g\big(t,x^0(t),u^0(t)\big)\big]=0.
\end{equation}
}
\end{sle}
\doc.\
   Докажем сначала, что это верно для моментов времени из построенной в теореме
   подпоследовательности $\tau'\subset\tau$.
   Как выше показано, для некоторых
   подпоследовательности $\tau'\subset\tau$ и последовательности
   пар $(\lambda^n,\psi^n)$,  для каждого $n\in\mm{N}$
    для почти всех $t<\tau''_n$ выполнено \rref{maxH}, то есть
          $$  \sup_{p\in
U(t)}\bigg(\psi^n(t)\big[f\big(t,x^0(t),p\big)-f\big(t,x^0(t),u^0(t)\big)\big]+
   \lambda^n\big[g\big(t,x^0(t),p\big)-g\big(t,x^0(t),u^0(t)\big)\big]\bigg)=0.$$
   Поскольку  $f\big(t,x^0(t),p\big)$
   ограничено в силу условия подлинейного роста, а   по построению
   имеет место
   $\psi^n(\tau'_n)=0$  (тогда, в частности, из \rref{dob} следует также $\lambda^n>0$),
   то,
   переходя к существенному пределу при $t\to\tau'_n-0$,
   имеем
   \rref{maxH_k} для данного $n\in\mm{N}$, то есть для всех
   моментов времени из
   подпоследовательности $\tau'.$

   Покажем, что \rref{maxH_k} выполнено для всех элементов из $\tau$,
   кроме, быть может, конечного их числа. Пусть не так, и
    \rref{maxH_k} не выполнено для счетного числа
   элементов. Тогда из них можно составить неограниченно возрастающую последовательность $\tau',$
   для каждого момента из которой, \rref{maxH_k} не имеет места.
   Поскольку управление $u^0$
   $\tau$-строго оптимально, то оно окажется и $\tau'$-строго
   оптимальным. В частности, для него найдется строго $\tau'$-исчезающее
   решение соотношений принципа максимума, а следовательно, по предыдущему абзацу,
   \rref{maxH_k} будет выполнено для каких-то элементов из
   $\tau'$, что противоречит выбору $\tau'.$
\bo

 Отметим также, что все показанное выше можно распространить на
 управляемые системы с негладкой правой частью.

 Отметим, что вообще говоря не утверждается, что $\psi^0(t)\to 0$ при $t\to\infty$,
 это может быть  не так даже вдоль всякой подпоследовательности
 последовательности $\tau$. Соответствующий пример приводится ниже.
\sect{Явная формула для сопряженной переменной}

 Воспользуемся тем, что система
 \rref{sys_psi} линейна,
  рассмотрим решение задачи Коши:
 \begin{equation*}
 \frac{d{A}(t)}{dt} =\frac{\partial f (t,x^0(t),u^0(t))}{\partial x}
  A(t),\quad A(0)=E.
\end{equation*}
  Введем теперь векторнозначную функцию $I$ времени
  правилом: для всех
  $T\in\mm{T}$
  \begin{equation*}
   I(T)\rav\int_0^T
   \frac{\partial g(t,x^0(t),u^0(t))}{\partial x}\, A(t)
  \,dt.
\end{equation*}
   Теперь всякое решение $\psi$ сопряженной системы \rref{sys_psi}
   для всех $t\in\mm{T}$ удовлетворяет формуле Коши
  \begin{equation}
   \label{4A}
   \psi(t)=(\psi(0)-\lambda I(t))A^{-1}(t).
\end{equation}

\begin{subequations}
\begin{pre}
\label{affp} {\it
     У всякого строгого $\tau$-исчезающего решения $(\psi^0,\lambda^0)$ всех соотношений принципа
     максимума \rref{sys_psi}--\rref{dob}
      для некоторой подпоследовательности $\tau'\subset\tau$,\\
    \begin{equation}
    \label{1p}
\text{либо}\quad  \lambda^0=1,\quad\psi^0(0)=  \lim_{n\to\infty}
I(\tau'^n),
\end{equation}
      и последовательность $I(\tau'^n)$ будет сходящейся и
      ограниченной,
%
    \begin{equation}
    \label{0p}
  \text{либо}\quad  \lambda^0=0,\quad  {\psi^0(0)}=
  \lim_{n\to\infty}\frac{I(\tau'^n)}{||I(\tau'^n)||},
  \end{equation}
   где данный предел существует, а последовательность $I(\tau'^n)$ неограничена.
\\
   Кроме того, для подпоследовательности $\tau'$ при всех $n\in\mm{N}$
   дополнительно выполнено  \rref{maxH_k}.
  }
\end{pre}
\end{subequations}
\doc.\

   Пусть некоторое решение $(\lambda^0,\psi^0)$  является
 строго  $\tau$-исчезающим решением; тогда $(\lambda^0,\psi^0)$ является пределом
   решений $(\lambda^n,\psi^n)$ системы
   \rref{sys_psi}, причем $\psi^n(\tau'^n)=0$
   для некоторой подпоследовательности $\tau'\subset\tau.$

   По формуле Коши  \rref{4A} для $t=\tau'^n$ из $\psi^n(\tau'^n)=0$ имеем
   $\psi^n(0)=\lambda^n I(\tau'^n),$
   отсюда
$$\psi^0(0)=\lim_{n\to\infty} \lambda^n I(\tau'^n).$$
   В случае  $\lambda^0\neq 0$ фактически показано равенство
   $\displaystyle \psi^0(0)= \lambda^0 \lim_{n\to\infty} I(\tau'^n).$
  Если же $\lambda^0=0,$  то, поскольку исчезающее решение нетривиально,
   получим $\psi^0(0)\neq 0$.
  Поделим полученное выражение на его же норму, имеем
  $$\frac{\psi(0)}{||\psi(0)||}=
  \lim_{n\to\infty}\frac{I(\tau'^n)}{||I(\tau'^n)||}.$$

\ \bo

 Полученные формулы можно несколько упростить.
\begin{zam}
\label{affss} {
\begin{subequations}
    Пусть существует конечный предел
    $$I_*\rav\lim_{n\to\infty}  I(\tau^n).$$
    Тогда выполнено \rref{1p}, и в которой
    $\psi^0$  задается формулой: для всех $T\in\mm{T}$
    \begin{equation}
    \label{1ss}
    \psi^0(T)=\big(I_*- I(T)\big)A^{-1}(T)=\\
     \Big(I_*- \int_{0}^T
   \frac{\partial g(t,x^0(t),u^0(t))}{\partial x}\, A(t)
  \,dt\Big)A^{-1}(T).
    \end{equation}

    Если же $I_*$ не зависит от выбора последовательности
    $\tau$, то есть существует конечный предел
        $\displaystyle\lim_{t\to\infty} I(t),$
        то последняя формула упрощается до
    \begin{equation}
    \label{1sss}
   \psi^0(T)=\int_{T}^\infty
   \frac{\partial g(t,x^0(t),u^0(t))}{\partial x}\, A(t)
  \,dt\, A^{-1}(T),
  \end{equation}
  где несобственный интеграл понимается в смысле Римана.
\end{subequations}
}
\end{zam}

   Отметим, что формула \rref{1sss} как явное выражение для
   сопряженной переменной  достаточно известна.
   Для линейных задач она была получена в работе \cite{aucl}.
   Для ряда
   стационарных задач как необходимое условие оптимальности она получена
   С.М.Асеевым, А.В.Кряжимским  в работах
   \cite{kr_as2,kr_as,kr_as3,kab};
   для некоторых нестационарных систем  её получили
   С.М.Асеев, V.M.Veliov (см.
   \cite{av,av_new}).
  Кроме вышеперечисленных работ,  условия, достаточные для применимости
  \rref{1sss}, рассматриваются также в работах автора
  \cite{mymy,mit,suz2012}. Там же есть условия для применимости формул \rref{1ss}, \rref{0p}.

   \medskip

   Приведем еще несколько отдельных замечаний для функционалов того
   или иного вида.

\begin{zam}
\label{affss900u} {
    Пусть функцию $g$ можно представить в виде
    $g(t,x,u)\rav g_1(t,x)-r(t)||u||^2,$ где функция $r$ положительна и ограничена на всяком компакте.
    Пусть кроме того, для почти всех
    $t\in\mm{T}$ выполнено $0\in U(t).$
    Тогда если $\tau$-строго оптимальное управление $u^0$
    существует, то выполнено $u(\tau_n-0)=0$ для всех $n\in\mm{N},$
    кроме, быть может, конечного их числа.
    }
\end{zam}
   Для доказательства замечания достаточно подставить такое $g$  в
  \rref{maxH_k}, и заметить, что только $0$ доставляет минимум $r(t)||u||^2$.

\begin{zam}
\label{affss121} {
    Пусть в условиях  замечания \ref{affss900u}  функцию $f$
    можно представить в виде $f(t,x,u)=f_1(t,x)+S(t)u$.
    Пусть кроме того, для каждого $n
    \in\mm{N}$  некоторая окрестность точки $0$   в  $\mm{X}$ содержится в  $S(t)U(t)$
    для почти всех   $t<n.$

    Тогда если $\tau$-строго оптимальное управление $u^0$
    существует, то       для всякого строго $\tau$-исчезающего решения
        $(\lambda^0,\psi^0)$,
       вместе с $u^0(\tau_n-0)$ зануляются
 $\psi^0(\tau_n)S(\tau_n-0)$ для всех $n\in\mm{N}$,
    кроме,
    быть может,
    конечного их числа. Более того, при этом $\lambda^0>0.$
    }
\end{zam}
 Действительно, для $\psi^0$ почти всюду выполнено  \rref{maxH}, подставляя такие $f,g$ имеем, что
    $u^0(t)$ почти всюду есть максимум (на $p\in U(t)$) функционала $\psi^0(t)S(t)p-\lambda^0||p||^2,$
 осталось заметить, что в силу предыдущего замечания  $u^0(\tau_n-0)=0$, а $0$~--- внутренняя точка образа $S\circ U$.

 \medskip

\begin{sle}
\label{affss900} {\it
    Пусть
    множество $$\Big\{\big(f(t,x^0(t),u),g(t,x^0(t),u)\big)\,\Big|\,u\in U(t)\Big\}$$ строго
    выпукло для $t$ из множества ненулевой меры $\ct{T}\subset\mm{T}$.

    Тогда для всякого  строго $\tau$-исчезающего решения $(\lambda^0,\psi^0)$
    соотношений принципа
    максимума \rref{sys_psi}--\rref{maxH}  выполнено также\\
    1) $\lambda^0>0,$\\
    2) $\psi^0(\tau'_n)=0$
     при всех $n\in\mm{N}$, кроме быть может, конечного их числа.
    }
\end{sle}
   \doc.\

      Действительно,
    для  строго $\tau$-исчезающего решения $(\lambda^0,\psi^0)$
    найдутся последовательность $(\lambda^n,\psi^n)$, сходящаяся к
    ней на всяком компакте, и соответствующая подпоследовательность
    $\tau'\subset\tau.$ Для каждого $n\in\mm{N}$, начиная с некоторого,
     найдется момент $T\in \ct{T}$ ($T<\tau'_n$), в который выполнено \rref{maxH} как для пары $(\lambda^0,\psi^0)$, так и для $(\lambda^n,\psi^n)$.
     Теперь в силу строгой выпуклости
    вектор $(\lambda^n,\psi^n(T))$ сонаправлен с $(\lambda^0,\psi^0(T))$, то есть они
   отличаются лишь на неотрицательный множитель.
    По формуле Коши \rref{4A}, отсюда, как значения $\psi^n(0),$ $\psi^0(0)$, так и сами
   функции $\psi^n,$ $\psi^0$, также отличаются на тот же множитель; в частности
   эти функции одновременно зануляются в момент времени $\tau'_n.$
   В силу выбора
     $n$ пункт 2) показан. В силу нетривиальности $\psi^0$
   как решения линейного уравнения, из 2) следует 1).
   \bo

 Отметим, что само по себе наличие формулы \rref{maxH_k} позволяет лишь написать дополнительную систему соотношений для выделения
 последовательности $\tau$. Последнее же утверждение указывает класс управляемых систем, в которых
 соотношения принципа максимума вместе со всеми краевыми условиями
 создают переопределенную систему соотношений, поскольку на каждый скалярный параметр $\tau_n$,
  в случае
 $dim \mm{X}>1$,  мы имеем больше одного соотношения.
 Именно этим, по-видимому, можно объяснить редкость строго оптимальных
 управлений.
\sect{Пример}

 Для демонстрации возможностей некоторых показанных выше утверждений
 рассмотрим модификацию примера Халкина
\cite{Halkin} (см. также
 \cite[Ex. 5.1]{Pickenhain},\, \cite[Ex.~1]{kr_as_t}):
 требуется  максимизировать
 функционал
   \begin{equation*}
   J[u](T)\rav\int_{[0,T]} (1-u)x\, dt
   \ \text{ при } \
   T\to\infty.
   \end{equation*}
   для управляемой системы
 \begin{eqnarray*}
\dot{x}=ux&\ &\text{для п.в.}\ t\in\mm{T};\\
x(0)=x_*>0,&\ & \\
u\in [\alpha,\beta]&\ & (\alpha\leqslant \beta).
   \end{eqnarray*}
 Будем искать в этом примере строго оптимальные управления.
 Каждому такому управлению $u^0$ соответствует некоторая
 неограниченно возрастающая последовательность
  $\tau$. В силу  \rref{maxH_k} можно считать, что $u^0(\tau_n-0)=\alpha$ для
  всех $n\in\mm{N}.$

 Легко видеть, что $A (T)=x^{0}(T)$ и $I(T)=
J[u^0](T).$  Переходя при необходимости от $\tau$ к ее
подпоследовательности, мы всегда можем считать выполненным вдоль
$\tau$ ровно один случай из трех:

   1) $J[u^0](\tau_n)\to -\infty$.  Теперь по формуле \rref{0p}
   имеем
  $\psi^0(0)=-1,$ $\lambda=0;$ подставляя найденное в гамильтониан получаем:
   $H=-u x,$ теперь из условия его максимума следует $u^0\equiv \alpha$.
     Проверкой убеждаемся, что
  для $J[u^0](\tau_n)\to -\infty$ необходимо $\alpha> 1$.

 2) $J[u^0](\tau_n)\to +\infty$; подобным образом получаем   $\psi^0(0)=-1,$ $\lambda=0;$
 $H=-u x,$ $u^0\equiv \beta,$ но это противоречит условию \rref{maxH_k} на
 выбор точек из $\tau.$

 3) $J[u^0](\tau_n)$ сходится к некоторому конечному $I_*$.
 Теперь можно воспользоваться \rref{1ss}, в частности считать $\lambda^0=1$.
  Введем функцию $R(t)\rav I_*-J[u^0](t)-x^0(t)$. Тогда
  $$H[t]=\big(I_*-J[u^0](t)\big)A^{-1}(t)u x^0(t)+(1-u)x^0(t)=
   \big(I_*-J[u^0](t)-x^0(t)\big)u+x^0(t)=R(t)u+x^0(t),$$ и $u^0(t)$
   в силу \rref{maxH}
   определяется знаком $R(t).$
   Поскольку  $\dot{R}(t)=-x(t)<0$, то $R$ монотонно убывает, а моментов переключения не больше
   одного.   Однако  $u^0$ принимает минимально возможное значение $\alpha$ в сколь
   угодно далекие  моменты времени, следовательно $R$ на больших
   временах отрицательно. Но тогда начиная с некоторого момента~$T$
   (или это момент переключения, или $T=0$) функция $R$ отрицательна,
   управление  $u^0$
   тождественно совпадает с $\alpha$,  а $x^0(t)=x^0(T)e^{\alpha(t-T)}$
   для всех $t>T$.

   С другой стороны,  последовательность чисел
   $I_*-J[u^0](\tau_n)$ сходится к  $0$, а последовательность
      чисел $R(\tau_n)$
   монотонно убывает, принимая начиная с некоторого, лишь отрицательные
   значения. Тогда их разность, последовательность чисел
   $x^0(\tau_n)=x^0(T)e^{\alpha(\tau_n-T)}$, не может сходиться к   $0$. Таким образом $\alpha >0.$

    В силу  фундаментальности последовательности чисел $I(\tau_n)=J[u^0](\tau_n)$  имеем
      для  их попарных разностей, чисел
     $\int_{[\tau_n,\tau_k]}(1-\alpha)e^{\alpha
     (t-\tau_n)}\,dt$ сходимость к  $0$ при больших $k,n$, что в силу $\alpha>0$ эквивалентно  $\alpha=1.$

      Но тогда, начиная с момента $T$, функция $J[u^0]$
     не изменяется, принимая значение $I_*$, следовательно
     $R(T)=-x(T)<0$. Теперь, если бы $T$ был моментом переключения
     управления $u^0$ (зависящего лишь от знака $R$), то $R(T)$ был бы равен
     $0$, следовательно переключения не было, тогда $T=0$, то есть  $u^0\equiv \alpha=1.$

     Заметим, что в случае $\alpha\geqslant 1$  управление $u^0\equiv
     \alpha$  является строго оптимальным.
     Действительно, если для пары управлений верно $u\geqslant v,$ то это верно и для их траекторий $x\geq
     y$, подставляя в целевой функционал, имеем $J[u]\leqslant J[v]$ в силу  $\alpha\geqslant 1$.
     Следовательно $u^0\equiv
     \alpha$, как  наименьшее из допустимых управлений, доставляет
     наибольший результат для каждого момента времени, в частности
     оно строго оптимально.

\medskip

   Отметим, что  для $\alpha>1$, хотя у всякого строго исчезающего решения
    сопряженная переменная $\psi^0$ является поточечным пределом зануляющихся
   во все более поздние моменты времени решений той же сопряженной системы, сама эта
   функция $\psi^0\equiv 1$ не стремится на бесконечности к нулю.
   Кроме того,  найденное  строго оптимальное решение
   является анормальным, в частности, формула \rref{1sss}, равно как и ее более общие модификации
   \rref{1p},\rref{1ss}, неприменима.
   Данный пример также показывает, что в замечании \ref{affss900}
   условие строгой выпуклости существенно.

\vspace{3ex}

\small

\makeatletter \@addtoreset{equation}{section}
\@addtoreset{footnote}{section}
\renewcommand{\section}{\@startsection{section}{1}{0pt}{1.3ex
plus 1ex minus 1ex}{1.3ex plus .1ex}{}}

{ 

\renewcommand{\refname}{{\rm\centerline{СПИСОК ЛИТЕРАТУРЫ}}}

\titleeng

\vspace{3ex}


\begin{center}
REFERENCES
\end{center}

1. Aseev S.M.,  Besov K.O.,  Kryazhimskii A.V. Infinite-horizon optimal
control problems in economics, {\it Russ. Math. Surv.,} 2012, Vol. 67, no. 2, pp. 195-253.

2. Aseev S.M.,  Kryazhimskii A.V. The Pontryagin maximum principle and
optimal economic growth problems, {\it Proc. Steklov Inst. Math.,} 2007, Vol. 257, pp. 1-
255.

3. Clarke F.H. Optimization and Nonsmooth Analysis, New York, J. Wiley, 1983.

4.  Pontryagin L.S.,  Boltyanskij V.G.,  Gamkrelidze R.V.,  Mishchenko E.F.
Mathematical Theory of Optimal Processes, Interscience Publishers, John Wiley and Sons,
New York, 1962.
5.  Khlopin D.V. On transversality condition in control problems for infinity horizon, {\it  Conference "Differential Equations and Optimal Control``, dedicated to the 90th anniversary of E. F. Mishchenko: Abstracts of Int. Conf.},
Steklov
Inst. Math., Moscow, 2008, pp. 144--146.

6.  Khlopin D.V.  On $\tau$-vanishing adjoint variable for infinity-horizon control problems,
 {\it  The International Conference on Differential Equations and Dynamical Systems: Abstracts of Int. Conf.},  Lomonosov Moscow State University, Moscow, 2012, pp. 173--174.

7.  { Aseev S.M., Kryazhimskii A.V.} The Pontryagin maximum principle
 and transversality conditions for a class of optimal control problems
 with infinite time horizons, {\it SIAM J. Control Optim.}, 2004, Vol. 43, pp. 1094–-1119.

8. { Aseev S.M.,  Kryazhimskii A.V.} Shadow prices in infinite-horizon
optimal control problems with dominating discounts, {\it Applied
Mathematics and Computation}, 2008, Vol.~204, no. 2,
pp. 519--531.

9. {  Aseev S.M.,  Kryazhimskii  A.V., Tarasyev A.M.}
The Pontryagin Maximum Principle and Transversality Conditions for an
Optimal Control Problem with Infinite Time Interval, {\it Proc. Steklov
Inst. Math.}, 2001, no. 233, pp. 64--80.

10. {  Aseev S.M.,  Veliov  V.M. }
Needle Variations in Infinite-Horizon
Optimal Control, {\it IIASA Interim Rept.
 IASA. Laxenburg, Research Report} 2012-04.~September 2012, 22~pp.

11. {  Aseev S.M.,  Veliov  V.M. } Maximum Principle for infinite-horizon
 optimal control problems with dominating discount, {\it Dynamics of Continuous, Discrete and Impulsive
Systems, Series B.}, 2012, V.~19, no. 1--2, pp. 43--63

12. { Aubin  J.P., Clarke F.H.} Shadow Prices and Duality for a Class
of Optimal Control Problems, {\it J. Control Optim.}, 1979,
Vol.17, pp. 567--586;

13.  { Bogucz D. } On the existence of a classical optimal solution and of an almost
  strongly optimal solution for an infinite-horizon control problem,
   {\it J. Optim. Theory Appl.},  2013, Vol. 156, no. 2.
  URL: http://link.springer.com/content/pdf/10.1007\%2Fs10957-012-0126-2
 DOI:  10.1007/s10957-012-0126-2

14. {  Carlson D.A. } Uniformly overtaking and
weakly overtaking optimal solutions in infinite-horizon optimal
control: when optimal solutions are agreeable {\it J. Optim. Theory
Appl.}, 1990, Vol.  64, no. 1, pp. 55--69.

15. {  Carlson D.A.,  Haurie A.B.,   Leizarowitz A.}, {\it Infinite Horizon Optimal
Control. Deterministic and Stochastic Systems.} Springer, Berlin,
1991.

16. Chakravarty S. The existence of an optimum savings program,  {\it
Econometrica}, 1962, Vol. 30, pp. 178–187

17. { Halkin H.} Necessary Conditions for Optimal Control Problems with
Infinite Horizons, {\it Econometrica}, 1974, Vol.~42,
pp. 267--272.

18. Khlopin D.V. Necessity of vanishing shadow price in infinite horizon control problems, { \it arXiv:1207.5358}
 URL: arxiv.org/pdf/1207.5358

19. {   Pickenhain S.} On adequate transversality conditions for
infinite horizon optimal control problems --- a famous example of
Halkin. In: Crespo Cuaresma, J.; Palokangas, T.; Tarasyev, A. (Eds.):
 {\it Dynamic Systems, Economic Growth, and the Environment.} Springer, 2010,  pp. 3--22

20. { Seierstad A. }
 Necessary conditions for nonsmooth, infinite-horizon optimal
 control problems, {\it J. Optim. Theory Appl.},  1999, Vol.~103, no.~1,
 pp. 201--230.

21. { Stern L.E.} Criteria of optimality in the infinite-time
optimal control problem, {\it J. Optim. Theory Appl.}, 1984, Vol.44, no 3, pp. 497-508.

22.  Weizacker C.C. Existence of optimal programs of
accumulation for an infinite time horizon, {\it Review of Economic Studies},  1965,
Vol. 32., pp.  85–104.

\received
\noindent
\contactinformationenglish

\end{document}